\documentclass{amsart}

\usepackage{amssymb,latexsym,amsfonts,amsmath}
\usepackage{graphicx}
\include{diagrams}
\usepackage{eso-pic}
\usepackage{epsfig}
\usepackage{graphicx}
\usepackage{color}
\usepackage{diagrams}
\usepackage{mathtools}
\usepackage{relsize}
\usepackage{tikz}
\usetikzlibrary{automata}
\usetikzlibrary{shapes}

\usepackage{amssymb,latexsym,amsfonts,amsmath}
\usepackage{graphicx}

\newtheorem{theorem}{Theorem}[section]

\theoremstyle{definition}
\newtheorem{definition}[theorem]{Definition}

\theoremstyle{remark}
\newtheorem{remark}[theorem]{Remark}
\numberwithin{equation}{section}

\newcommand{\R}{{\mathbb{R}}}

\newcommand{\N}{{\mathbb{N}}}

\begin{document}

\begin{abstract}
Finite-state models of control systems were proposed by several researchers as a convenient mechanism to synthesize controllers enforcing complex specifications. Most techniques for the construction of such symbolic models have two main drawbacks: either they can only be applied to restrictive classes of systems, or they require the exact computation of reachable sets. In this paper, we propose a new abstraction technique that is applicable to any smooth control system as long as we are only interested in its behavior in a compact set. Moreover, the exact computation of reachable sets is not required. The effectiveness of the proposed results is illustrated by synthesizing a controller to steer a vehicle.
\end{abstract}

\title[Symbolic models for nonlinear control systems without stability assumptions]{Symbolic models for nonlinear control systems\\ without stability assumptions}
\thanks{This work has been partially supported by the National Science Foundation award 0717188, 0820061, European Commission under STREP project HYCON$^2$, and by the Center of Excellence for Research DEWS, University of LÕAquila, Italy.}

\author[majid zamani]{Majid Zamani$^1$} 
\author[giordano pola]{Giordano Pola$^2$} 
\author[Manuel Mazo Jr.]{Manuel Mazo Jr.$^3$} 
\author[paulo tabuada]{Paulo Tabuada$^1$}
\address{$^1$Department of Electrical Engineering\\
University of California at Los Angeles,
Los Angeles, CA 90095}
\email{\{zamani, tabuada\}@ee.ucla.edu}
\urladdr{http://www.ee.ucla.edu/~zamani}
\urladdr{http://www.ee.ucla.edu/~tabuada}
\address{$^2$Department of Electrical and Information Engineering, Center of Excellence DEWS, University of LÕAquila, 
Poggio di Roio, 67040 LÕAquila, Italy}
\email{giordano.pola@univaq.it}
\urladdr{http://www.diel.univaq.it/people/pola}
\address{$^3$INCAS$^3$, Dr. Nassaulaan 9, 9401 HJ Assen, The Netherlands and the Faculty of Mathematics and Natural Sciences, ITM, University of Groningen,
Groningen, 9747AG, The Netherlands}
\email{M.Mazo@rug.nl}
\urladdr{http://www.rug.nl/staff/m.mazo/index}

\maketitle

\section{Introduction}
In the past years several different abstraction techniques have been developed to assist in the synthesis of controllers enforcing complex specifications. This paper is concerned with symbolic abstractions resulting from replacing aggregates or collections of states of a control system by symbols. When a symbolic abstraction with a finite number of states or symbols is available, the synthesis of the controllers can be reduced to a fixed-point computation over the finite-state abstraction \cite{paulo}. Moreover, by leveraging computational tools developed for discrete-event systems \cite{DEDSBook,CassandrasBook} and games on automata \cite{InfGames,ComputeGames,AVW03}, one can synthesize controllers satisfying specifications difficult to enforce with conventional control design methods. Examples of such specification classes include logic specifications expressed in linear temporal logic or automata on infinite strings.

The quest for symbolic abstractions has a long history including results on timed automata \cite{alur}, rectangular hybrid automata \cite{puri}, and o-minimal hybrid systems \cite{lafferriere, brihaye}. Early results for classes of control systems were based on dynamical consistency properties \cite{caines}, natural invariants of the control system \cite{koutsoukos}, $l$-complete approximations \cite{moor}, and quantized inputs and states \cite{forstner,BMP02}. Recent results include work on piecewise-affine and multi-affine systems \cite{habets,BH06}, set-oriented discretization approach for discrete-time nonlinear optimal control problem \cite{junge1}, abstractions based on an elegant use of convexity of reachable sets for sufficiently small time \cite{gunther}, and the use of incremental input-to-state stability \cite{pola,pola1,pola2,girard2}.

Our results improve upon most of the existing techniques in two directions: i) by being applicable to larger classes of control systems; ii) by not requiring the exact computation of reachable sets which is a hard task in general. In the first direction, our technique improves upon the results in \cite{BMP02,habets, BH06} by being applicable to systems not restricted to non-holonomic chained-form, piecewise-affine, and multi-affine systems, respectively, and upon the results in \cite{pola, pola1, pola2, girard2} by not requiring any stability assumption. In the second direction, our technique improves upon the results in \cite{moor,forstner} by not requiring the exact computation of reachable sets. The results in \cite{junge1} offer a discretization tailored to optimal control while our discretization is independent of the control objective. In \cite{gunther} a different abstraction technique is proposed that is also applicable to a wide class of control systems and does not require the exact computation of reachable sets. Such technique provides tight over-approximations of reachable sets based on convexity but requires small sampling times. Other efficient techniques are available in the literature for computing over-approximations of reachable sets. For example, \cite{junge,dellnitz,patrick} provide tight over-approximations of reachable sets, not necessarily convex, at the cost of a higher computational complexity than \cite{gunther}. In contrast to \cite{gunther,patrick}, our technique imposes no restrictions on the choice of the sampling time but provides less tight over-approximations of the set of reachable states.

In this paper, we show that symbolic models exist if the control systems satisfy an \textit{incremental forward completeness} assumption which is an incremental version of forward completeness. The main contribution of this paper is to establish that: 

\emph{For every nonlinear control system satisfying the incremental forward completeness assumption, one can construct a symbolic model that is alternatingly approximately simulated \cite{pola1} by the control system and that approximately simulates \cite{girard} the control system.} Although these results are of theoretical nature, we also provide a simple way of constructing symbolic models which can be improved by using tighter over-approximations of reachable sets such as those described in~\cite{gunther,junge,dellnitz}.

We illustrate the results presented in this paper through a simple example in which a vehicle is requested to reach a target set while avoiding a number of obstacles.
  
\section{Control Systems and Incremental Forward Completeness\label{sec2}}
\subsection{Notation} 
The identity map on a set $A$ is denoted by $1_{A}$. If $A$ is a subset of $B$ we denote by \mbox{$\imath
_{A}:A\hookrightarrow B$} or simply by $\imath $ the natural inclusion map taking any $a \in A$ to \mbox{$\imath (a) = a \in B$}. The symbols $\mathbb{N}$, $\mathbb{Z}$, $\mathbb{R}$, $\mathbb{R}^+$ and $\mathbb{R}_0^+$ denote the set of natural, integer, real, positive, and nonnegative real numbers, respectively. Given a vector \mbox{$x\in\mathbb{R}^{n}$}, we denote by $x_{i}$ the $i$--th element of $x$, and by $\Vert x\Vert$ the infinity norm of $x$. Given a matrix $M\in\mathbb{R}^{n\times{m}}$, we denote  by $\Vert{M}\Vert$ the infinity norm of $M$. The closed ball centered at $x\in{\mathbb{R}}^{n}$ with radius $\varepsilon$ is defined by \mbox{$\mathcal{B}_{\varepsilon}(x)=\{y\in{\mathbb{R}}^{n}\,|\,\Vert x-y\Vert\leq\varepsilon\}$}. For any set \mbox{$A\subseteq\R^n$} of the form $A=\bigcup_{j=1}^MA_j$ for some $M\in\N$, where $A_j=\prod_{i=1}^n [c_i^j,d_i^j]\subseteq \R^n$ with $c^j_i<d^j_i$ and positive constant $\eta\leq\widehat\eta$, where $\widehat{\eta}=\min_{j=1,\cdots,M}\eta_{A_j}$ and \mbox{$\eta_{A_j}=\min\{|d_1^j-c_1^j|,\cdots,|d_n^j-c_n^j|\}$}, define \mbox{$[A]_{\eta}=\{a\in A\,\,|\,\,a_{i}=k_{i}\eta,k_{i}\in\mathbb{Z},i=1,\cdots,n\}$}. The set $[A]_{\eta}$ will be used as an approximation of the set A with precision $\eta$. Note that $[A]_{\eta}\neq\varnothing$ for any $\eta\leq\widehat\eta$. Geometrically, for any $\eta\in{\mathbb{R}^+}$ and $\lambda\geq\eta$ the collection of sets \mbox{$\{\mathcal{B}_{\lambda}(p)\}_{p\in[A]_{\eta}}$} is a covering of $A$, i.e. \mbox{$A\subseteq\bigcup_{p\in[A]_{\eta}}\mathcal{B}_{\lambda}(p)$}. By defining \mbox{$[\R^n]_{\eta}=\{a\in \R^n\,\,|\,\,a_{i}=k_{i}\eta,k_{i}\in\mathbb{Z},i=1,\cdots,n\}$}, the set \mbox{$\bigcup_{p\in[\R^n]_{\eta}}\mathcal{B}_{\lambda}(p)$} is a covering of $\R^n$ for any $\eta\in\R^+$ and $\lambda\geq\eta/2$. Given a measurable function \mbox{$f:\mathbb{R}_{0}^{+}\rightarrow\mathbb{R}^n$}, the (essential) supremum ($\sup$ norm) of $f$ is denoted by $\Vert f\Vert_{\infty}$; we recall that $\Vert{f}\Vert_\infty=(\text{ess})\sup\left\{\Vert{f(t)}\Vert,t\geq0\right\}$. A continuous function \mbox{$\gamma:\mathbb{R}_{0}^{+}\rightarrow\mathbb{R}_{0}^{+}$}, is said to belong to class $\mathcal{K}$ if it is strictly increasing and \mbox{$\gamma(0)=0$}; function $\gamma$ is said to belong to class $\mathcal{K}_{\infty}$ if \mbox{$\gamma\in\mathcal{K}$} and $\gamma(r)\rightarrow\infty$ as $r\rightarrow\infty$. 

\subsection{Control Systems\label{II.B}}

The class of control systems that we consider in this paper is formalized in
the following definition.

\begin{definition}
\label{Def_control_sys}A \textit{control system} $\Sigma$ is a quadruple $\Sigma=(\mathbb{R}^{n},\mathsf{U},\mathcal{U},f)$, where:
\begin{itemize}
\item $\mathbb{R}^{n}$ is the state space;
\item $\mathsf{U}\subseteq\R^m$ is the input set; 
\item $\mathcal{U}$ is a subset of all piecewise continuous functions of time from intervals of the form \mbox{$]a,b[\subseteq\mathbb{R}$} to $\mathsf{U}$ with $a<0$ and $b>0$; 
\item \mbox{$f:\mathbb{R}^{n}\times \mathsf{U}\rightarrow\mathbb{R}^{n}$} is a continuous map
satisfying the following Lipschitz assumption: for every compact set
\mbox{$Q\subset\mathbb{R}^{n}$}, there exists a constant $L\in\mathbb{R}^+$ such that for all $x,y\in Q$ and all $u\in \mathsf{U}$, we have \mbox{$\Vert
f(x,u)-f(y,u)\Vert\leq L\Vert x-y\Vert$}.
\end{itemize}
\end{definition}

A curve \mbox{$\xi:]a,b[\rightarrow\mathbb{R}^{n}$} is said to be a
\textit{trajectory} of $\Sigma$ if there exists $\upsilon\in\mathcal{U}$
satisfying \mbox{$\dot{\xi}(t)=f\left(\xi(t),\upsilon(t)\right)$}, for almost all $t\in$ $]a,b[$. We also write $\xi_{x\upsilon}(\tau)$ to denote the point reached at time $\tau$
under the input $\upsilon$ from initial condition $x=\xi_{x\upsilon}(0)$; this point is
uniquely determined, since the assumptions on $f$ ensure existence and
uniqueness of trajectories \cite{sontag1}. Although we have defined trajectories over open domains, we shall refer to
trajectories \mbox{${\xi_{x\upsilon}:}[0,\tau]\rightarrow\mathbb{R}^{n}$} and input curves $\upsilon:[0,\tau[\to \mathsf{U}$, with the understanding of the
existence of a trajectory \mbox{${\xi}_{x\upsilon'}^{\prime}:]a,b[\rightarrow\mathbb{R}^{n}$} and input curve $\upsilon':]a,b[\to \mathsf{U}$ such that \mbox{${\xi}_{x\upsilon}={\xi}_{x\upsilon'}^{\prime}|_{[0,\tau]}$} and $\upsilon=\upsilon'\vert_{[0,\tau[}$. Note that by continuity of $\xi$, $\xi_{x\upsilon}(\tau)$ is uniquely defined as the left limit of $\xi_{x\upsilon}(t)$ with $t\rightarrow\tau$.

A control system $\Sigma$ is said to be forward complete if every trajectory is defined on an interval of the form $]a,\infty[$. Sufficient and necessary conditions for a system to be forward complete can be found in \cite{sontag}. 

\subsection{Incremental forward completeness}
The results presented in this paper require a certain property that we introduce in this section. 
\begin{definition}
\label{dFC}
A control system $\Sigma$ is incrementally forward complete ($\delta$-FC) if it is forward complete and there exist continuos functions \mbox{$\beta: \mathbb{R}_0^+\times\mathbb{R}_0^+\rightarrow\mathbb{R}_0^+$} and \mbox{$\gamma: \mathbb{R}_0^+\times\mathbb{R}_0^+\rightarrow\mathbb{R}_0^+$} such that for every $s\in\mathbb{R}^+$, the functions $\beta(\cdot,s)$ and $\gamma(\cdot,s)$ belong to class $\mathcal{K}_{\infty}$, and for any $x,x'\in{\mathbb{R}^n}$, any $\tau\in\R^+$, and any \mbox{$\upsilon, \upsilon'\in{\mathcal{U}}$}, where \mbox{$\upsilon, \upsilon':[0,\tau[\rightarrow{\mathsf{U}}$}, the following condition is satisfied for all $t\in[0,\tau]$:
\begin{equation}
\left\Vert\xi_{x\upsilon}(t)-\xi_{x'\upsilon'}(t)\right\Vert\leq\beta(\Vert{x}-x'\Vert,t)+\gamma(\left\Vert{\upsilon}-\upsilon'\right\Vert_{\infty},t). \label{uns_cond}%
\end{equation}
\end{definition}
Incremental forward completeness requires the distance between two arbitrary trajectories to be bounded by the sum of two terms capturing the mismatch between the initial conditions and the mismatch between the inputs as shown in (\ref{uns_cond}). 

\begin{remark}
We note that $\delta$-FC implies uniform continuity of the map $\phi_t:\R^n\times\mathcal{U}\to \R^n$ defined by $\phi_t(x,\upsilon)=\xi_{x\upsilon}(t)$ for any fixed $t\in \R_0^+$. Here, uniform continuity is understood with respect to the topology induced by the infinity norm on $\R^n$, the $\sup$ norm on $\mathcal{U}$, and the product topology on $\R^n\times \mathcal{U}$.
\end{remark}

Note that a linear control system:
\begin{equation}
\dot{\xi}=A\xi+B\upsilon,~~\xi(t)\in\mathbb{R}^n,~\upsilon(t)\in \mathsf{U}\subseteq\mathbb{R}^m, \nonumber
\end{equation}
is $\delta$-FC and the functions $\beta$ and $\gamma$ can be chosen as:
\begin{equation}
\beta(r,t)=\left\Vert e^{At}\right\Vert{r};~~\gamma(r,t)=\left(\int_0^{t}\left\Vert{e}^{As}B\right\Vert ds\right)r,
\end{equation}
where $\Vert{e}^{At}\Vert$ denotes the infinity norm of $e^{At}$.

The notion of $\delta$-FC can be described in terms of Lyapunov-like functions.
We start by introducing the following definition which was inspired by the notion of incremental input-to-state stability ($\delta$-ISS) Lyapunov function presented in \cite{angeli}.
\begin{definition}
\label{delta_FC}
Consider a control system $\Sigma$ and a smooth function $V:\mathbb{R}^n\times\mathbb{R}^n\rightarrow\mathbb{R}_0^+$. Function $V$ is called a $\delta$-FC Lyapunov function for $\Sigma$, if there exist $\mathcal{K}_{\infty}$ functions $\underline{\alpha}$, $\overline{\alpha}$, $\sigma$, and $\kappa\in\mathbb{R}$ such that:
\begin{itemize}
\item[(i)] for any $x,x'\in\mathbb{R}^n$, $\underline{\alpha}(\Vert{x}-x'\Vert)\leq{V}(x,x')\leq\overline{\alpha}(\Vert{x}-x'\Vert)$;
\item[(ii)] for any $x,x'\in\mathbb{R}^n$ and for any $u,u'\in\mathsf{U}$, $\frac{\partial{V}}{\partial{x}}f(x,u)+\frac{\partial{V}}{\partial{x'}}f(x',u')\leq \kappa V(x,x')+\sigma(\Vert{u}-u'\Vert)$.
\end{itemize}
\end{definition}
The following theorem describes $\delta$-FC in terms of the existence of a $\delta$-FC Lyapunov function.
\begin{theorem}
\label{theorem4}
A control system $\Sigma=(\R^n,\mathsf{U},\mathcal{U},f)$ is $\delta$-FC if it admits a $\delta$-FC Lyapunov function. Moreover, the functions $\beta$ and $\gamma$ in (\ref{uns_cond}) are given by:
 \begin{eqnarray}
 \beta(r,t)=\underline{\alpha}^{-1}\left(2e^{\kappa{t}}\overline{\alpha}(r)\right),~~\gamma(r,t)=\underline{\alpha}^{-1}\left(2\frac{e^{\kappa{t}}-1}{\kappa}\sigma(r)\right). \label{beta&gamma}
 \end{eqnarray}
\end{theorem}
The proof of the preceding result is reported in \cite{majid} and was inspired by the work in \cite{sontag}.

\section{Symbolic Models and Approximate Equivalence Notions}\label{symbolic}
\subsection{Systems and control systems}
We use systems to describe both control systems as well as their symbolic models. A more detailed exposition of the notion of system that we now introduce can be found in~\cite{paulo}.
\begin{definition}\cite{paulo}
A system $S$ is a quintuple $S=(X,U,\longrightarrow,Y,H)$ consisting of:

\begin{itemize}
\item A set of states $X$;
\item A set of inputs $U$;
\item A transition relation $\longrightarrow\subseteq X\times U\times X$;
\item An output set $Y$;
\item An output function $H:X\rightarrow Y$.
\end{itemize}
\end{definition}
System $S$ is said to be:

\begin{itemize}
\item \textit{metric}, if the output set $Y$ is equipped with a metric
$\mathbf{d}:Y\times Y\rightarrow\mathbb{R}_{0}^{+}$;
\item \textit{countable}, if $X$ is a countable set;
\item \textit{finite}, if $X$ is a finite set.
\end{itemize}

A transition \mbox{$(x,u,x')\in\longrightarrow$} is denoted by $x\rTo^ux'$. For a transition $x\rTo^ux'$, state $x'$ is called a \mbox{$u$-successor}, or simply successor, of state $x$. We denote by $\mathbf{Post}_{u}(x)$ the set of \mbox{$u$-successors} of a state $x$ and by $U(x)$ the set of inputs $u\in{U}$ for which $\mathbf{Post}_{u}(x)$ is nonempty. We shall abuse the notation and denote by $\mathbf{Post}_{u}(Z)$ the set $\mathbf{Post}_{u}(Z)=\bigcup_{x\in{Z}}\mathbf{Post}_{u}(x)$. A system is deterministic if for any state $x\in{X}$ and any input $u$, there exists at most one \mbox{$u$-successor} (there may be none). A system is called nondeterministic if it is not deterministic. Hence, for a nondeterministic system it is possible for a state to have two (or possibly more) distinct $u$-successors.

\begin{definition}\cite{paulo}
For a system $S=(X,U,\longrightarrow,Y,H)$ and given any state $x_0\in{X}$, a finite state run generated from $x_0$ is a finite sequence of transitions:
$$
x_0\rTo^{u_0}x_1\rTo^{u_1}x_2\rTo^{u_2}\cdots\rTo^{u_{n-2}}x_{n-1}\rTo^{u_{n-1}}x_n,
$$
 such that $x_i\rTo^{u_i}x_{i+1}$ for all $0\leq i<n$. In some cases, a finite state run can be extended to an infinite state run. 
 
 An infinite state run generated from $x_0$ is an infinite sequence:
$$x_0\rTo^{u_0}x_1\rTo^{u_1}x_2\rTo^{u_2}x_3\rTo^{u_3}\cdots$$
such that $x_i\rTo^{u_i}x_{i+1}$ for all $i\in\N_0$.
 \end{definition}

\subsection{System relations}
We start by recalling approximate simulation relations, introduced in \cite{girard}, that are useful when analyzing or synthesizing controllers for deterministic systems. 
\begin{definition}
\label{ASR}Let \mbox{$S_{a}=(X_{a},U_{a},\rTo_{a},Y_a,H_{a})$} and
\mbox{$S_{b}=(X_{b},U_{b},\rTo_{b},Y_b,H_{b})$} be metric systems
with the same output sets $Y_a=Y_b$ and metric $\mathbf{d}$, and consider a precision $\varepsilon\in\mathbb{R}^{+}$. A
relation \mbox{$R\subseteq X_{a}\times X_{b}$} is said to be an $\varepsilon
$-approximate simulation relation from $S_{a}$ to $S_{b}$, if the following three conditions are satisfied:

\begin{itemize}
\item[(i)] for every $x_{a}\in{X_a}$, there exists $x_{b}\in{X_b}$ with $(x_{a},x_{b})\in{R}$;

\item[(ii)] for every $(x_{a},x_{b})\in R$ we have \mbox{$\mathbf{d}(H_{a}(x_{a}),H_{b}(x_{b}))\leq\varepsilon$};

\item[(iii)] for every $(x_{a},x_{b})\in R$ we have that \mbox{$x_{a}\rTo_{a}^{u_a}x'_{a}$ in $S_a$} implies the existence of \mbox{$x_{b}\rTo_{b}^{u_b}x'_{b}$} in $S_b$ satisfying $(x'_{a},x'_{b})\in R$.
\end{itemize}  

System $S_{a}$ is \mbox{$\varepsilon$-approximately} simulated by $S_{b}$ or $S_b$ \mbox{$\varepsilon$-approximately} simulates $S_a$, denoted by \mbox{$S_{a}\preceq_{\mathcal{S}}^{\varepsilon}S_{b}$}, if there exists
an \mbox{$\varepsilon$-approximate} simulation relation from $S_{a}$ to $S_{b}$.
\end{definition}

For nondeterministic systems we need to consider relationships that explicitly capture the adversarial nature of nondeterminism. The notion of alternating approximate simulation relation is shown in \cite{pola1} to be appropriate to this regard. 

\begin{definition}
\label{AASR} Let $S_a$ and $S_b$ be metric systems with the same output sets $Y_a=Y_b$ and metric $\mathbf{d}$, and consider a precision $\varepsilon\in\mathbb{R}^+$. A relation \mbox{$R\subseteq X_a\times X_b$} is said to be an $\varepsilon$-approximate alternating simulation relation from $S_{a}$ to $S_{b}$ if conditions (i), (ii) in Definition \ref{ASR} and the following condition are satisfied:
\begin{itemize}
\item[(iii)] for every $(x_a,x_b)\in R$ and for every $u_a\in U_a(x_a)$ there exists $u_b\in U_b(x_b)$ such that for every $x'_b \in \mathbf{Post}_{u_b}(x_b)$ there exists \mbox{$x'_a \in \mathbf{Post}_{u_a}(x_a)$} satisfying \mbox{$(x'_{b},x'_{a})\in R$}.
\end{itemize}
\end{definition}
System $S_a$ is alternatingly \mbox{$\varepsilon$-approximately} simulated by $S_b$ or $S_b$ alternatingly \mbox{$\varepsilon$-approximately} simulates $S_a$, denoted by \mbox{$S_{a}\preceq_{\mathcal{A}\mathcal{S}}^{\varepsilon}S_{b}$}, if there exists an alternating \mbox{$\varepsilon$-approximate} simulation relation from $S_a$ to $S_b$. 

It is readily seen from the above definitions that the notions of approximate simulation and of alternating approximate simulation coincide when the systems involved are deterministic. 

The importance of the preceding notions lies in enabling the transfer of controllers designed for a symbolic model to controllers acting on the original control system. More details about these notions and how the refinement of controllers can be performed are reported in \cite{paulo}.

\section{Symbolic Models for $\delta$-FC Control Systems}\label{existence}
This section contains the main contribution of the paper. We show that the time discretization of a $\delta$-FC control system, suitably restricted to a compact set, admits a finite abstraction.

The results in this section rely on additional assumptions on $\mathsf{U}$ and $\mathcal{U}$ that we now describe. Such assumptions are not required for the definitions and results in Sections \ref{sec2} and \ref{symbolic}. We restrict attention to control systems \mbox{$\Sigma=(\mathbb{R}^{n},\mathsf{U},\mathcal{U},f)$} with input sets $\mathsf{U}$ of the form $\mathsf{U}=\bigcup_{j=1}^J\mathsf{U}_j$ for some $J\in\N$, where $\mathsf{U}_j=\prod_{i=1}^m [a_i^j,b_i^j]\subseteq \R^m$ with $a^j_i<b^j_i$. For such input sets we define the constant $\widehat{\mu}=\min_{j=1,\cdots,J}\mu_{{\mathsf{U}_j}}$ where $\mu_{{\mathsf{U}_j}}=\min\{|b_1^j-a_1^j|,\cdots,|b_m^j-a_m^j|\}$. We further restrict attention to sampled-data control systems, where input curves belong to $\mathcal{U}_\tau$ containing only constant curves of duration $\tau\in\R^+$, i.e.
\begin{equation}
 \nonumber
 \mathcal{U}_\tau=\{\upsilon:[0,\tau[\to \mathsf{U}\,\,\vert\,\,\upsilon(t)=\upsilon(0),t\in[0,\tau[\}.
 \end{equation}
 
Given a sampling time $\tau\in \R^+$ and a control system \mbox{$\Sigma=(\mathbb{R}^{n},\mathsf{U},\mathcal{U}_\tau,f)$}, consider the system $S_{\tau}(\Sigma)=(X_{\tau},U_{\tau},\rTo_{\tau},Y_{\tau},H_{\tau})$ consisting of:
\begin{itemize}
\item $X_{\tau}=\mathbb{R}^{n}$;
\item $U_{\tau}=\mathcal{U}_\tau$;
\item $x_{\tau}\rTo_{\tau}^{\upsilon_{\tau}}x'_{\tau}$ if there exists a trajectory
\mbox{$\xi_{x_{\tau}\upsilon_{\tau}}:[0,\tau]\rightarrow\mathbb{R}^{n}$} of $\Sigma$ satisfying \mbox{$\xi_{x_{\tau}\upsilon_{\tau}}(\tau)=x'_{\tau}$};
\item $Y_{\tau}=\mathbb{R}^{n}$;
\item $H_{\tau}=1_{\mathbb{R}^{n}}$.
\end{itemize}

The above system can be thought of as the time discretization of the control system $\Sigma$. Indeed, a finite state run
\[
x_{0}\rTo^{\upsilon_{1}}_{\tau}x_{1}\rTo^{\upsilon_{2}}_{\tau} \,...\, \rTo^{\upsilon_{N}}_{\tau} x_{N}
\]
of $S_{\tau}(\Sigma)$ captures the state evolution of the control system $\Sigma$ at times $t=0,\tau,\cdots,N\tau$. The state run starts from the initial condition $x_{0}$, with control input $\upsilon$, obtained by the concatenation of control inputs $\upsilon_{i}$ \big(i.e. $\upsilon(t)=\upsilon_{i}(0)$ for any $t\in [(i-1)\tau,i\,\tau[$\big), for $i=1,\cdots,N$.

We consider a $\delta$-FC control system \mbox{$\Sigma=(\mathbb{R}^{n},\mathsf{U},\mathcal{U}_\tau,f)$}, and a quadruple $\mathsf{q}=(\tau,\eta,\mu,\theta)$ of quantization parameters, where $\tau\in\mathbb{R}^{+}$ is the sampling time, $\eta\in\mathbb{R}^{+}$ is the state space quantization, $\mu\in\mathbb{R}^{+}$ is the input set quantization, and $\theta\in\mathbb{R}^{+}$ is a design parameter. Define the system:
\begin{equation}
S_{\mathsf{q}}(\Sigma)=(X_{\mathsf{q}},U_{\mathsf{q}},\rTo_{\mathsf{q}},Y_{\mathsf{q}},H_{\mathsf{q}}),
\label{T2}%
\end{equation}
consisting of:
\begin{itemize}
\item $X_{\mathsf{q}}=[\mathbb{R}^{n}]_{\eta};$
\item $U_{\mathsf{q}}=[\mathsf{U}]_{\mu};$
\item $x_{\mathsf{q}}\rTo_{\mathsf{q}}^{u_{\mathsf{q}}}x'_{\mathsf{q}}$ if $\Vert\xi_{x_{\mathsf{q}}u_{\mathsf{q}}}(\tau)-x'_{\mathsf{q}}\Vert\leq\beta(\theta,\tau)+\gamma\left(\mu,\tau\right)+\eta$;
\item $Y_{\mathsf{q}}=\mathbb{R}^{n}$;
\item $H_{\mathsf{q}}=\imath : X_{\mathsf{q}}\hookrightarrow Y_{\mathsf{q}}$,
\end{itemize}
where $\beta$ and $\gamma$ are the functions appearing in (\ref{uns_cond}). In the definition of the transition relation, and in the remainder of the paper, we abuse notation by identifying $u_\mathsf{q}$ with the constant input curve with domain $[0,\tau[$ and value $u_{\mathsf{q}}$.

The transition relation of $S_\mathsf{q}(\Sigma)$ is well defined in the sense that for every $x_\mathsf{q}\in X_\mathsf{q}$ and every $u_\mathsf{q}\in U_\mathsf{q}$ there always exists $x_\mathsf{q}'\in X_\mathsf{q}$ such that $x_\mathsf{q}\rTo_{\mathsf{q}}^{u_\mathsf{q}}x_\mathsf{q}'$. This can be seen by noting that by definition of $X_\mathsf{q}$, for any $x\in \R^n$ there always exists a state $x_\mathsf{q}'\in X_\mathsf{q}$ such that $\Vert x-x_\mathsf{q}'\Vert\leq\eta$. Hence, for $x=\xi_{x_\mathsf{q}u_\mathsf{q}}(\tau)$ there always exists a state $x_\mathsf{q}'\in X_\mathsf{q}$ satisfying \mbox{$\Vert\xi_{x_{\mathsf{q}}u_{\mathsf{q}}}(\tau)-x'_{\mathsf{q}}\Vert\leq\eta\leq\beta(\theta,\tau)+\gamma\left(\mu,\tau\right)+\eta$}.

We can now state the main result of the paper which relates $\delta$-FC to existence of symbolic models. 

\begin{theorem}
Let $\Sigma=(\R^n,\mathsf{U},\mathcal{U}_\tau,f)$ be a $\delta$-FC control system. For any desired precision $\varepsilon\in\mathbb{R}^+$, and any quadruple $\mathsf{q}=(\tau,\eta,\mu,\theta)$ of quantization parameters satisfying $\mu\leq\widehat\mu$ and $\eta\leq\varepsilon\leq\theta$, we have:
\begin{equation}
S_{\mathsf{q}}(\Sigma)\preceq_{\mathcal{A}\mathcal{S}}^{\varepsilon}S_{\tau}(\Sigma)\preceq_{\mathcal{S}}^{\varepsilon}S_{\mathsf{q}}(\Sigma).
\end{equation}
\label{theorem1}
\end{theorem}

\begin{proof}
We start by proving \mbox{$S_{\tau}(\Sigma)\preceq_{\mathcal{S}}^{\varepsilon}S_{\mathsf{q}}(\Sigma)$}.
Consider the relation $R\subseteq X_{\tau}\times X_{\mathsf{q}}$ defined by \mbox{$(x_{\tau},x_{\mathsf{q}})\in R$}
if and only if \mbox{$\Vert H_{\tau}(x_{\tau})-H_{\mathsf{q}}(x_{\mathsf{q}})\Vert=\Vert x_{\tau}-x_{\mathsf{q}}\Vert\leq\varepsilon$}. Since \mbox{$X_{\tau}\subseteq\bigcup_{p\in[\mathbb{R}^n]_{\eta}}\mathcal{B}_{\eta}(p)$}, for every $x_{\tau}\in{X_{\tau}}$ there exists \mbox{$x_{\mathsf{q}}\in{X}_{\mathsf{q}}$} such that: 
\begin{equation}
\Vert{x_{\tau}}-x_{\mathsf{q}}\Vert\leq\eta\leq\varepsilon. 
\end{equation}
Hence, \mbox{$(x_{\tau},x_{\mathsf{q}})\in{R}$} and condition (i) in Definition \ref{ASR} is satisfied. Now consider any \mbox{$(x_{\tau},x_{\mathsf{q}})\in R$}. Condition (ii) in Definition \ref{ASR} is satisfied
by the definition of $R$. Let us now show that condition (iii) in Definition
\ref{ASR} holds.

Consider any \mbox{$\upsilon_{\tau}\in {U}_{\tau}$}. Choose an input \mbox{$u_{\mathsf{q}}\in U_{\mathsf{q}}$} satisfying:
\begin{equation}
\Vert \upsilon_{\tau}-u_{\mathsf{q}}\Vert_{\infty}=\Vert \upsilon_{\tau}(0)-u_{\mathsf{q}}(0)\Vert\leq\mu.\label{b01}%
\end{equation}
Note that the existence of such $u_\mathsf{q}$ is guaranteed by the special shape of $\mathsf{U}$, described in the beginning of this section, and by the inequality $\mu\le\widehat{\mu}$ which guarantees that $\mathsf{U}\subseteq\bigcup_{p\in[\mathsf{U}]_{\mu}}\mathcal{B}_{{\mu}}(p)$. Consider the unique transition \mbox{$x_{\tau}\rTo_{\tau}^{\upsilon_{\tau}}x'_{\tau}=\xi_{x_{\tau}\upsilon_{\tau}}(\tau)$} in $S_{\tau}(\Sigma)$. It follows from the $\delta$-FC assumption that the distance between $x'_{\tau}$ and  \mbox{$\xi_{x_{\mathsf{q}}u_{\mathsf{q}}}(\tau)$} is bounded as:
\begin{equation}
\Vert x'_{\tau}-\xi_{x_{\mathsf{q}}u_{\mathsf{q}}}(\tau)\Vert\leq\beta(\varepsilon,\tau)+\gamma\left(\mu,\tau\right).\label{b02}%
\end{equation}
Since \mbox{$X_{\tau}\subseteq\bigcup_{p\in[\mathbb{R}^n]_{\eta}}\mathcal{B}_{\eta}(p)$}, there exists \mbox{$x'_{\mathsf{q}}\in{X}_{\mathsf{q}}$} such that: 
\begin{equation}
\Vert{x'_{\tau}}-x'_{\mathsf{q}}\Vert\leq\eta. \label{b04}
\end{equation}
Using the inequalities $\varepsilon\leq\theta$, (\ref{b02}), and (\ref{b04}), we obtain:
\begin{align*}
\Vert{\xi_{x_{\mathsf{q}}u_{\mathsf{q}}}(\tau)}-x'_{\mathsf{q}}\Vert&\leq\Vert{\xi_{x_{\mathsf{q}}u_{\mathsf{q}}}(\tau)}-{x'_{\tau}}\Vert+\Vert{x'_{\tau}}-x'_{\mathsf{q}}\Vert\\ \notag
&\leq\beta(\varepsilon,\tau)+\gamma\left(\mu,\tau\right)+\eta\leq\beta(\theta,\tau)+\gamma\left(\mu,\tau\right)+\eta, \notag
\end{align*}
which, by the definition of $S_\mathsf{q}(\Sigma)$, implies the existence of \mbox{$x_{\mathsf{q}}\rTo_{\mathsf{q}}^{u_{\mathsf{q}}}x'_{\mathsf{q}}$ in $S_{\mathsf{q}}(\Sigma)$}. Therefore, from inequality (\ref{b04}) and since $\eta\leq\varepsilon$, we conclude \mbox{$(x'_{\tau},x'_{\mathsf{q}})\in{R}$} and condition (iii) in Definition \ref{ASR} holds. 

Now we prove \mbox{$S_{\mathsf{q}}(\Sigma)\preceq_{\mathcal{AS}}^{\varepsilon}S_{\tau}(\Sigma)$}. Consider the relation \mbox{$R\subseteq X_{\tau}\times X_{\mathsf{q}}$}, defined in the first part of the proof. For every \mbox{$x_{\mathsf{q}}\in{X}_{\mathsf{q}}$}, by choosing \mbox{$x_{\tau}=x_{\mathsf{q}}$}, we have \mbox{$(x_{\tau,}x_{\mathsf{q}})\in{R}$} and condition (i) in Definition \ref{AASR} is satisfied. Now consider any \mbox{$(x_{\tau},x_{\mathsf{q}})\in R$}. Condition (ii) in Definition \ref{AASR} is satisfied
by the definition of $R$. Let us now show that condition (iii) in Definition
\ref{AASR} holds. 
Consider any \mbox{$u_{\mathsf{q}}\in U_{\mathsf{q}}$}. Choose the input $\upsilon_\tau=u_\mathsf{q}$ and consider the unique $x'_{\tau}=\xi_{x_{\tau}\upsilon_{\tau}}(\tau)\in\mathbf{Post}_{\upsilon_{\tau}}(x_\tau)$ in $S_{\tau}(\Sigma)$.
From the $\delta$-FC assumption, the distance between $x'_{\tau}$ and \mbox{$\xi_{x_{\mathsf{q}}u_{\mathsf{q}}}(\tau)$} is bounded as:
\begin{equation}
\Vert x'_{\tau}-\xi_{x_{\mathsf{q}}u_{\mathsf{q}}}(\tau)\Vert\leq\beta(\varepsilon,\tau).\label{b03}%
\end{equation}

Since \mbox{$X_{\tau}\subseteq\bigcup_{p\in[\mathbb{R}^n]_{\eta}}\mathcal{B}_{\eta}(p)$}, there exists \mbox{$x'_{\mathsf{q}}\in{X}_{\mathsf{q}}$} such that: 
\begin{equation}
\Vert{x'_{\tau}}-x'_{\mathsf{q}}\Vert\leq\eta. \label{b05}
\end{equation}
Using the inequalities, $\varepsilon\leq\theta$, (\ref{b03}), and (\ref{b05}), we obtain:
\begin{align*}
\Vert{\xi_{x_{\mathsf{q}}u_{\mathsf{q}}}(\tau)}-x'_{\mathsf{q}}\Vert\leq\Vert{\xi_{x_{\mathsf{q}}u_{\mathsf{q}}}(\tau)}-{x'_{\tau}}\Vert+\Vert{x'_{\tau}}-x'_{\mathsf{q}}\Vert\leq\beta(\varepsilon,\tau)+\eta\leq\beta(\theta,\tau)+\gamma\left(\mu,\tau\right)+\eta, \notag
\end{align*}
which, by definition of $S_\mathsf{q}(\Sigma)$, implies the existence of \mbox{$x_{\mathsf{q}}\rTo_{\mathsf{q}}^{u_{\mathsf{q}}}x'_{\mathsf{q}}$ in $S_{\mathsf{q}}(\Sigma)$}. Therefore, from inequality (\ref{b05}) and since $\eta\leq\varepsilon$, we can conclude that \mbox{$(x'_{\tau},x'_{\mathsf{q}})\in{R}$} and condition (iii) in Definition \ref{ASR} holds. 
\end{proof}

\begin{remark}\label{remark5}
Whenever $\mathcal{U}_\tau$ only contains finite number of curves, the function $\gamma$ is not required to construct $S_{\mathsf{q}}(\Sigma)$. This can be seen by noting that we can use all the elements in $\mathcal{U}_\tau$ when constructing $S_{\mathsf{q}}(\Sigma)$ thus eliminating the approximation error on input curves, represented by the term $\gamma(\mu,\tau)$ in the definition of $\rTo_{\mathsf{q}}$.
\end{remark}

\begin{remark}
\label{remark1}
The transition relation defined in (\ref{T2}) can also be written as:
\begin{equation}
x_{\mathsf{q}}\rTo_{\mathsf{q}}^{u_{\mathsf{q}}}x'_{\mathsf{q}}~~~\text{if}~~~\mathcal{B}_{\eta}(x'_\mathsf{q})\cap\mathcal{B}_{\beta(\theta,\tau)+\gamma\left(\mu,\tau\right)}(\xi_{x_{\mathsf{q}}u_{\mathsf{q}}}(\tau))\neq\varnothing.
\end{equation} 
This shows that we place a transition from $x_\mathsf{q}$ to any point $x_\mathsf{q}'$ for which the ball $\mathcal{B}_{\eta}(x_\mathsf{q}')$ intersects the over-approximation of $\mathbf{Post}_{u_\mathsf{q}}\left(\mathcal{B}_\varepsilon(x_\mathsf{q})\right)$ in $S_\tau(\Sigma)$ given by $\mathcal{B}_{\beta(\theta,\tau)+\gamma\left(\mu,\tau\right)}(\xi_{x_{\mathsf{q}}u_{\mathsf{q}}}(\tau))$. It is not difficult to see that the conclusion of Theorem \ref{theorem1} remains valid if we use any other over-approximation of the set \mbox{$\mathbf{Post}_{u_\mathsf{q}}(\mathcal{B}_{\varepsilon}(x_\mathsf{q}))$} in $S_\tau(\Sigma)$.
\end{remark}

The symbolic model $S_\mathsf{q}(\Sigma)$ has a countably infinite set of states. In order to construct a finite symbolic model we note that in practical applications the physical variables are restricted to a compact set. Velocities, temperatures, pressures, and other physical quantities cannot become arbitrarily large without violating the operational envelop defined by the control problem being
solved. By making use of this fact, we can directly compute a finite abstraction $S_{\mathsf{q}D}(\Sigma)$ of $S_\tau(\Sigma)$ capturing the behavior of $S_\tau(\Sigma)$ within a given set $D$ of the form \mbox{$D=\bigcup_{j=1}^MD_j$} for some $M\in\N$, where $D_j=\prod_{i=1}^n [c_i^j,d_i^j]\subseteq \R^n$ with $c^j_i<d^j_i$, describing the valid range for the physical variables. By having the extra condition $\eta\leq\widehat\eta$, where $\widehat{\eta}=\min_{j=1,\cdots,M}\eta_{D_j}$ where \mbox{$\eta_{D_j}=\min\{|d_1^j-c_1^j|,\cdots,|d_n^j-c_n^j|\}$}, we define the system $S_{\mathsf{q}D}(\Sigma)=(X_{\mathsf{q}D},U_{\mathsf{q}D},\rTo_{\mathsf{q}D},Y_{\mathsf{q}D},H_{\mathsf{q}D})$, where $U_{\mathsf{q}D}=U_\mathsf{q}$, $Y_{\mathsf{q}D}=Y_\mathsf{q}$, and $H_{\mathsf{q}D}=H_\mathsf{q}$ and
\begin{itemize}
\item $X_{\mathsf{q}D}=[D]_\eta$;
\item $x_{\mathsf{q}D}\rTo_{\mathsf{q}D}^{u_{\mathsf{q}D}}x'_{\mathsf{q}D}$ if $\Vert\xi_{x_{\mathsf{q}D}u_{\mathsf{q}D}}(\tau)-x'_{\mathsf{q}D}\Vert\leq\beta(\theta,\tau)+\gamma\left(\mu,\tau\right)+\eta$ and any $x'_\mathsf{q}\in\mathbf{Post}_{u_{\mathsf{q}D}}(x_{\mathsf{q}D})$ in $S_\mathsf{q}(\Sigma)$ belongs to $X_{\mathsf{q}D}$;
\end{itemize}
Note that $S_{\mathsf{q}D}(\Sigma)$ is a finite system because $D$ is a compact set. Moreover, the relation $R\subseteq X_{\mathsf{q}D}\times X_{\mathsf{q}}$ defined by $(x_{\mathsf{q}D},x_{\mathsf{q}})\in R$ if \mbox{$x_{\mathsf{q}D}=x_{\mathsf{q}}$} is a $0$-approximate alternating simulation relation from $S_{\mathsf{q}D}(\Sigma)$ to $S_{\mathsf{q}}(\Sigma)$.  By combining $S_{\mathsf{q}D}(\Sigma)\preceq_{\mathcal{AS}}^0 S_\mathsf{q}(\Sigma)$ with $S_\mathsf{q}(\Sigma)\preceq_\mathcal{AS}^\varepsilon S_\tau(\Sigma)$ we conclude\footnote{It is shown in~\cite{paulo} that the composition of two alternating simulation relations is still an alternating simulation relation.} $S_{\mathsf{q}D}(\Sigma)\preceq_\mathcal{AS}^\varepsilon S_\tau(\Sigma)$. Hence, any controller synthesized for the finite model $S_{\mathsf{q}D}(\Sigma)$ can be refined to a controller enforcing the same specification on $S_\tau(\Sigma)$. Detailed information on how to construct refinements can be found in \cite{paulo}.

\section{Example}
We illustrate the results of the paper on a vehicle. We borrowed this example from \cite{astrom}. In this model, the motion of the front and rear pairs of wheels are approximated by a single front wheel and a single rear wheel.
We consider the following model for the vehicle:
\begin{equation}\label{vehicle}
    \Sigma:\left\{
                \begin{array}{ll}
                  \dot{x}=v_0\frac{\cos(\alpha+\theta)}{\cos(\alpha)},\\ 
                 \dot{y}=v_0\frac{\sin(\alpha+\theta)}{\cos(\alpha)},\\ 
                 \dot{\theta}=\frac{v_0}{b}\tan(\delta),
                \end{array}
                \right.
\end{equation}
where $\alpha=\arctan\left(\frac{a\tan(\delta)}{b}\right)$. The position of the vehicle is given by the pair $(x,y)$, and the orientation of the vehicle is given by $\theta$. The pair $(v_0,\delta)$ are the control inputs, expressing the velocity of the rear wheel and the steering angle, respectively. It is readily seen that $\Sigma$ is not incrementally input-to-state stable \cite{angeli}. Hence, the results in \cite{pola, pola1} cannot be applied to this system. We assume that $a=0.5$, $b=1$, $(v_0,\delta)\in{\mathsf{U}}=[-1,~1]\times[-1,~1]$ and that the control inputs are piecewise constant. Since control inputs are piecewise constant of duration $\tau$, it can be readily checked that for any $t\in[0,\tau]$, we get:
\begin{eqnarray}
\nonumber
x(t)&=&\frac{b}{\cos(\alpha)\tan(\delta)}\Big[\sin\left(\alpha+\frac{v_0}{b}\tan(\delta)t+\theta(0)\right)-\sin(\alpha+\theta(0))\Big]+x(0),\\\notag
y(t)&=&\frac{b}{\cos(\alpha)\tan(\delta)}\Big[\cos\left(\alpha+\frac{v_0}{b}\tan(\delta)t+\theta(0)\right)-\cos(\alpha+\theta(0))\Big]+y(0),\\\notag
\theta(t)&=&\frac{v_0}{b}\tan(\delta)t+\theta(0),
\end{eqnarray}
if $\tan(\delta)\neq0$, and
\begin{eqnarray}
\nonumber
x(t)=v_0\cos(\theta(0))t+x(0),~~y(t)=v_0\sin(\theta(0))t+y(0),~~\theta(t)=\theta(0),
\end{eqnarray}
if $\tan(\delta)=0$. It can be verified that for the given $\mathsf{U}$, the function $\beta$ is given by $\beta(r,t)=(1+1.267t)r$. Here we are assuming that $\mathcal{U}_\tau$ is finite and contains curves taking values in $[\mathsf{U}]_{0.3}$. Hence, as explained in Remark \ref{remark5}, the function $\gamma$ is not required to construct the abstraction.
 
\begin{figure}
\centering
\includegraphics[width=10cm]{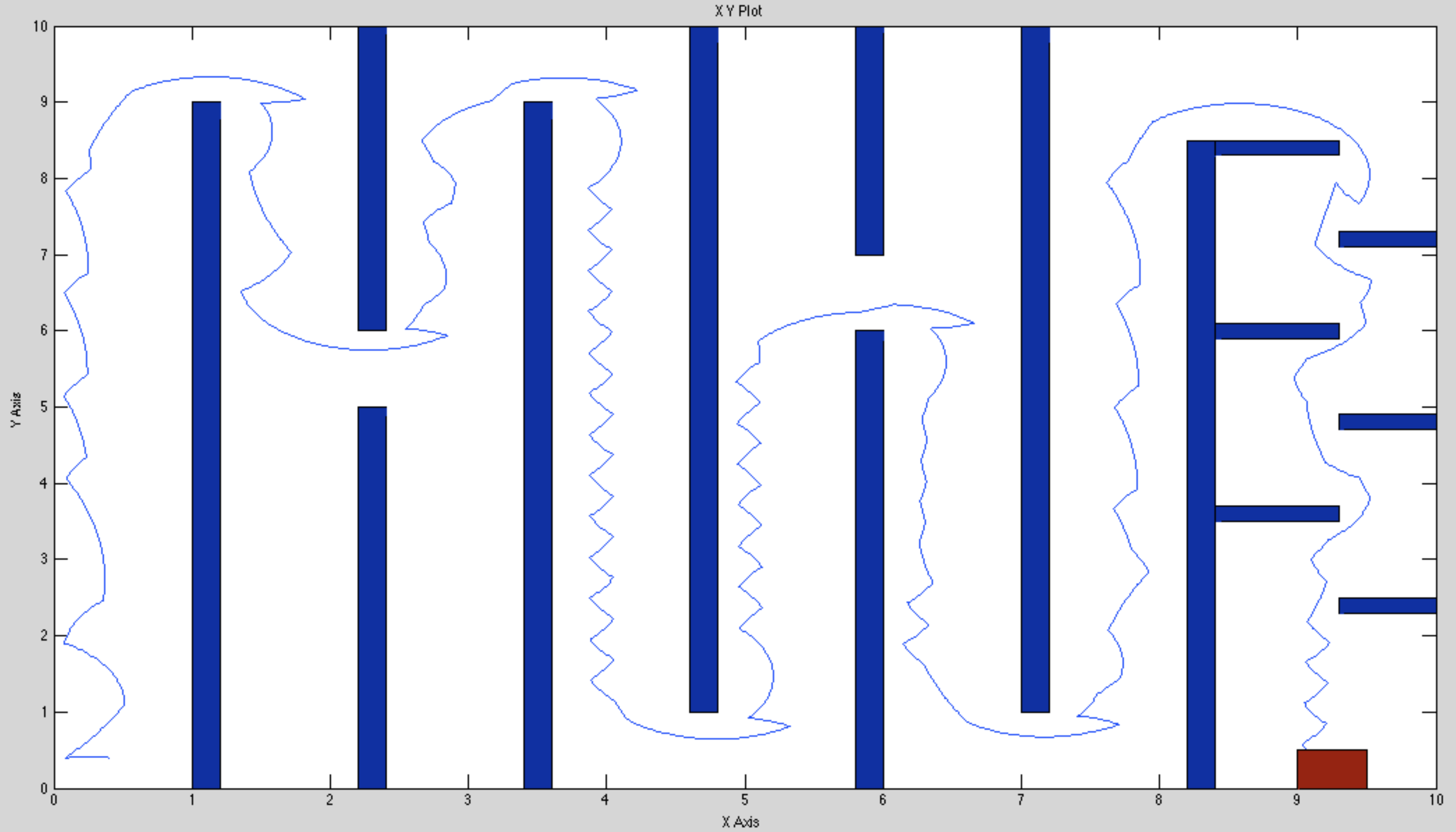}
\caption {Evolution of the vehicle with initial condition \mbox{$(0.4,~0.4,~0)$}.}\label{spec1}
\end{figure}

We work on the subset $D=[0,~10]\times[0,~10]\times[-\pi,~\pi]$ of the state space of $\Sigma$. Our objective is to design a controller navigating the vehicle to reach the target set $W=[9,~9.5]\times[0,~0.5]$, indicated with a red box in Figure \ref{spec1}, while avoiding the obstacles, indicated as blue boxes in Figure \ref{spec1}, and remain indefinitely inside $W$.

For a precision $\varepsilon=0.2$, we construct a symbolic model $S_{\mathsf{q}D}(\Sigma)$ by choosing $\theta=0.2$, $\eta=0.2$, and $\tau=0.3$ so that assumptions of Theorem \ref{theorem1} are satisfied. The computation of the abstraction $S_{\mathsf{q}D}(\Sigma)$ was performed using the tool\footnote{\textsf{Pessoa} can be freely downloaded from \texttt{http://www.cyphylab.ee.ucla.edu/pessoa}.} \textsf{Pessoa}~\cite{pessoa}. A controller enforcing the specification has been found by using standard algorithms from game theory, see e.g. \cite{paulo}. 

\begin{figure}
\centering
\includegraphics[width=10cm]{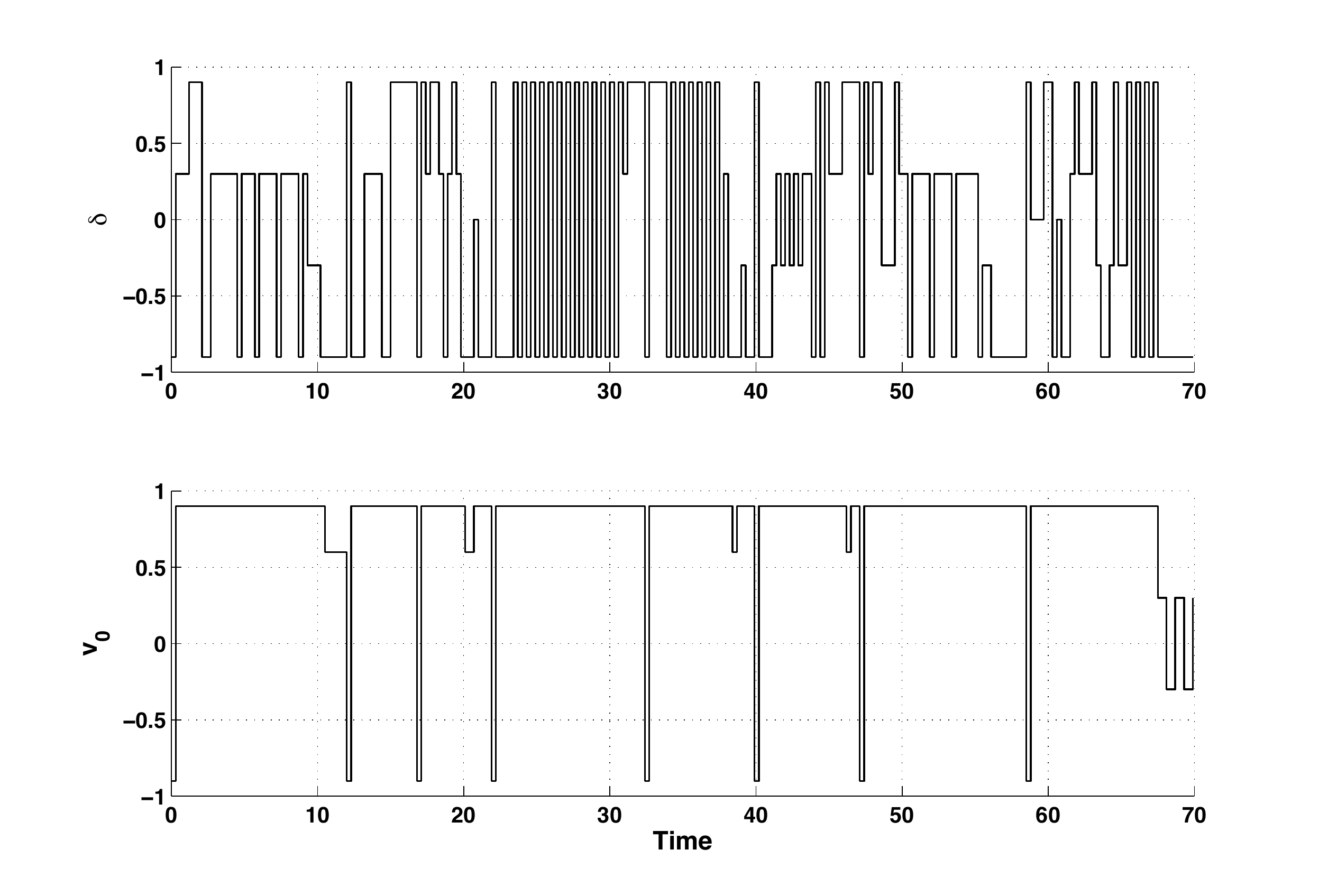}
\caption {Evolution of the input signals.}\label{input}
\end{figure}

In Figure \ref{spec1}, we show the closed-loop trajectory stemming from the initial condition \mbox{$(0.4,~0.4,~0)$}. It is readily seen that the specification are satisfied.
In Figure \ref{input}, we show the evolution of input signals.

\section{Discussion}
In this paper we showed that any smooth control system, suitably restricted to a compact subset of states, admits a finite symbolic model. The proposed symbolic model can be used to synthesize controllers enforcing complex specifications given in several different formalisms such as temporal logics or automata on infinite strings. The synthesis of such controllers is well understood and can be performed using simple fixed-point computations as described in~\cite{paulo}. The current limitation of this design methodology is the size of the computed abstractions. The authors are currently investigating several different techniques to address this limitation such as integrating the design of controllers with the construction of symbolic models~\cite{pola3}. Efforts by other researchers include the use of non-uniform quantization~\cite{tazaki}.

\bibliographystyle{alpha}
\bibliography{reference}
\end{document}